\documentclass{article}
\usepackage{amsfonts}
\usepackage[letterpaper,body={14cm,22cm}, mag=1000]{geometry}
\usepackage{amssymb}
\usepackage{secdot}
\usepackage{setspace}
\usepackage{titlesec}
\titleformat{\section}[runin]{\bfseries\filcenter}{\thesection}{1em}{}

\renewcommand{\thesection}{\arabic{section}}
\title{\large \bf Class-preserving automorphisms of some finite $p$-groups}
\author{\small \bf Mahak Sharma \\
\small \em  Department of Applied Sciences and Humanities\\
\small \em Baddi University of Emerging Sciences and Technology, Baddi, India\\
\\
\small \bf Deepak Gumber\footnote{Research supported by University Grants Commission of India under the Research Award Scheme}\\
\small \em  School of Mathematics and Computer Applications\\
\small \em  Thapar University, Patiala - 147 004, India\\}

\date{}
\newtheorem{thm}{Theorem}[section]
\newtheorem{lm}[thm]{Lemma}
\newtheorem{pp}[thm]{Proposition}

\newtheorem{cor}[thm]{Corollary}

\begin{document}

\maketitle

\begin{abstract}
\noindent {\bf Abstract.} Let $G$ be a finite $p$-group of order $p^5$, where $p$ is a prime. We give necessary and sufficient conditions on $G$ such that $G$ has a non-inner class-preserving automorphism. As a consequence, we give short and alternate proofs of results of section 5 of Yadav [Proc. Indian Acad. Sci. (Math. Sci.) {\bf 118} (2008), 1-11] and Theorem 4.2 of Kalra and Gumber [Indian J. Pure Appl. Math. {\bf 44} (2013), 711-725].
\end{abstract}

\vspace{1ex}


\vspace{1ex}


\section{Introduction}
Let $G$ be a finite group. An automorphism $\alpha$ of $G$ is called a class-preserving  automorphism if for each element $x \in G $,
there exists an element $g_{x}\in G$ such that $\alpha(x)=g_{x}^{-1}xg_{x}$; and is called an inner automorphism if for all $ x\in G$,
there exists a fix element $g \in G $ such that $ \alpha(x)=g^{-1}xg$. The group $\mathrm{Inn}(G)$ of all inner automorphisms of $G$ is a normal subgroup of the group $\mathrm{Aut}_{c}(G)$ of all class-preserving automorphisms of $G$. We denote the group $ \mathrm{Aut}_{c}(G)/\mathrm{Inn}(G)$  of all class-preserving outer
automorphisms by $\mathrm{Out}_{c}(G)$. The interest in class-preserving automorphisms took place
when Burnside \cite[p. 463]{Bur1} asked the question: Does there exist any finite group $G$ such that $G$ has a non-inner class-preserving automorphism? Burnside \cite{Bur2} himself gave an affirmative
answer to his question by constructing a group $G$ of order $p^{6}$ for which $\mathrm{Out}_{c}(G)$ was non-trivial. For more details about this
problem, one can see the survey article by Yadav \cite{Yad2011}.

In this note, we are especially interested in class-preserving automorphisms of $p$-groups
of order upto $p^{5}$. Using the known classifications and presentations of finite extra-special $p$-groups and of groups of order $p^4$, Kumar and Vermani \cite {Yad2000, Yad2001} proved that if  $G$ is an extra-special $p$-group (in particular if $|G| = p^{3}$) or if $|G|=p^{4}$, then $\mathrm{ Out}_{c}(G) = 1 $.  It follows from \cite{Jam} that all finite $p$-groups of order $p^5$, where $p$ is an odd prime, are partitioned into ten isoclinism families; and from \cite{HalSen} that all finite 2-groups of order $2^5$ are partitioned into eight isoclinism families. Yadav \cite{Yad2008} proved that if $G$ and $H$ are two finite non-abelian isoclinic groups, then $\mathrm{Aut}_c(G)\simeq\mathrm{Aut}_c(H)$. He then showed, by picking up one group from each isoclinism family, that $\mathrm{Out}_c(G)\ne 1$ for the groups $\Phi_7(1^5)$ and $\Phi_{10}(1^5)$ from seventh and tenth family in \cite{Jam}, and hence concluded that if $G$ is a finite $p$-group of order $p^5$, where $p$ is an odd prime, then  $\mathrm{Out}_c(G)\ne 1$ if and only if $G$ is isoclinic to a group either in the seventh or in the tenth family. Recently, Kalra and Gumber \cite[Theorem 4.2]{KalGum}, using the classifications, isoclinism families and presentations given by Hall and Senior \cite{HalSen} and Sag and Wamsley \cite{SagWam}, have shown that if $G$ is a finite 2-group of order $2^5$, then   $\mathrm{Out}_c(G)=1$ except for the forty fourth and forty fifth groups from the sixth family in \cite{HalSen}. 

Hertweck \cite[Proposition 14.4]{Her} proved that if $G$ is a finite group having an abelian normal subgroup $A$ with cyclic quotient $G/A$, then class-preserving automorphisms of $G$ are inner automorphisms. Yadav \cite[Corollary 3.6]{Yad2008} proved that if $G$ is a finite $p$-group of nilpotence class $2$ such that $G'$ is cyclic, then $\mathrm{Out}_c(G)=1$. We shall use these results quite frequently in our proofs. Observe that if $G$ is an extra-special $p$-group, then $G'$ is cyclic; and if $|G|=p^4$, then $G$ has a maximal abelian subgroup. It follows that for these groups $G$, $\mathrm{Out}_c(G)=1$. The present note is the result of an effort to find necessary and sufficient conditions, without using the available classifications and presentations, on a finite $p$-group $G$ for which $\mathrm{Out}_c(G)$ is non-trivial. In section 2, we prove our main theorem, Theorem 2.3, which gives the necessary and sufficient conditions on a finite $p$-group $G$ of order $p^{5}$ for which $\mathrm{Out}_{c}(G) \neq 1 $. As a consequence, we obtain short and alternate  proofs of the results of Yadav \cite[Section 5]{Yad2008} and Kalra and Gumber \cite[Theorem 4.2]{KalGum}.

 An automorphism $ \alpha $ of a group $G$ is called a central automorphism if it induces the identity
 automorphism on $G/Z(G)$; or equivalently, $x^{-1}\alpha(x)\in Z(G)$ for all $x\in G$. By $\mathrm{Aut}_z(G)$ we denote the group of all central automorphisms of $G$. For $x \in G$, $x^{G}$ denotes the conjugacy class of $x$ in $G$.
By $G'$ and $\Phi(G)$, we respectively denote the commutator and the Frattini subgroup of G. A cyclic group of order $m$ is denoted as $C_{m}$ and $\gamma_3(G)$ denotes the third term of the lower central series of $G$. The symbol $cl(G)$ denotes the nilpotence class of $G$ and by $d(G)$ we denote the smallest cardinality of a generating set of $G$. If $H$ is a non-trivial proper normal subgroup of $G$, then $(G,H)$ is called a Camina pair if and only if $ H \subseteq [x,G]$ for all $x \in G-H$, where
$[x,G]$ = $\{[x, g]|g \in G\}$.

\section{Proof of Theorem.} We start with the following crucial technical lemma. The lemma can be of independent interest also.

\begin{lm}
Let $G$ be any group such that $\mathrm{Out}_{c}(G/Z(G))$ is trivial. Then $\mathrm{Aut}_{c}(G)=(\mathrm{Aut}_c(G)\cap \mathrm{Aut}_{z}(G))\mathrm{Inn}(G)$. In addition, if $G$ is finite, then
$$|\mathrm{Aut}_{c}(G)|=\frac{|\mathrm{Aut}_{c}(G)\cap \mathrm{Aut}_{z}(G)|\;|\mathrm{Inn}(G)|}{|Z(\mathrm{Inn}(G))|}\cdot$$
\end{lm}
{\bf Proof.}
 If $\alpha \in \mathrm{Aut}_{c}(G)$, then it induces a class-preserving automorphism, say $\overline{\alpha}$, on $G/Z(G)$ given by $\overline{\alpha}(xZ(G))=\alpha(x)Z(G)$ for all $x\in G$. Since $\mathrm{Out}_{c}(G/Z(G))=1$, $\alpha(x)Z(G)=
a^{-1}xaZ(G)$ for a fix $a\in G$. Therefore, for each $x \in G $, there exists an element $z_{x} \in Z(G)$ such that $\alpha (x)  =  a^{-1}xaz_{x} = a^{-1}(xz_{x})a$. Define a map $\beta:G \rightarrow G$ by $\beta(x)=a \alpha(x) a^{-1}$ for all $ x \in G $. 
It is easy to see that $\beta \in \mathrm{Aut}_c(G)\cap \mathrm{Aut}_z(G)$ and $ \alpha = i_{a}\beta $, where $i_{a}$ denotes the inner automorphism of $G$ given by conjugation with $a$. It thus follows that $\mathrm{Aut}_{c}(G)= (\mathrm{Aut}_c(G) \cap \mathrm{ Aut}_{z}(G)) \mathrm{Inn}(G)$.
\hfill $\Box$

\begin{pp}Let $G$ be a finite non-abelian $p$-group of order $p^{5}$ such that $|Z(G)|\ge p^2$. Then $\mathrm{Out}_{c}(G)=1 $.
\end{pp}
{\bf Proof.} If $|Z(G)|\ge p^3$, then $G$ has a maximal abelian subgroup and hence $\mathrm{Out}_{c}(G) = 1 $. We therefore suppose that $|Z(G)|$= $p^{2}$. Then $cl(G)$ is either 2 or 3. First suppose that $cl(G)$ = 2.
Then $G' \leq Z(G)$. If $ G'$ is cyclic, then $\mathrm{Out}_{c}(G) = 1 $. Therefore
suppose that $G' = Z(G) \simeq C_{p} \times C_{p}$. Then $\mathrm{exp}(G/Z(G)) = \mathrm{exp}(G') = p$. It follows that $Z(G) = G'=\Phi(G)$ and hence $d(G)=3$. If $G$ = $\langle a,b,c \rangle$, then   $G'$= $\langle [a,b], [a,c], [b,c] \rangle$.
 We can assume that $[b,c] = [a,b]^{m}[a,c]^{n}$ for some $m,n$, $0 \leq m,n\leq p-1$.  Set $u:= ba^{-n}$ and $v:= ca^{m}$. Then 
 $$[u,v] = [ba^{-n},ca^{m}] = [b,c][b,a]^{m}[a,c]^{-n}=1.$$ It follows that $ \langle u,v, G'\rangle $ is a maximal 
abelian subgroup of $G$ and hence $\mathrm{Out}_{c}(G)=1$.

 Now suppose that $cl(G)$ = 3. Then $ G'$ is abelian of order $p^{2}$ or $p^{3}$. First assume that $|G'|=p^{2}$. Since $G/C_{G}(G')$ is isomorphic to a subgroup of $\mathrm{Aut}( G')$, $[G :C_{G}(G')] \leq p$ and hence $ |C_{G}(G')| = p^{4}$. The subgroup $G'Z(G)$ is of order $p^3$ and is contained in $C_G(G')$. It follows that $C_G(G')=\langle a,G'Z(G)\rangle$, where $a\in C_G(G')-G'Z(G)$, is a maximal abelian subgroup of $G$ and hence $\mathrm{Out}_{c}(G)=1$. Next assume that $|G'|=p^{3}$. Then $G'=\Phi(G)$  and $d(G)=2$. Let $G=\langle a,b\rangle$ and let $w:=[a,b],u:=[a,w]$ and $v:=[b,w]$. Then $G'=\langle w,\gamma_3(G) \rangle $  and $\gamma_3(G)=\langle u,v\rangle$. Since $[a,b]^{p}\equiv [a,b^{p}] \equiv 1 \pmod{ \gamma_{3}(G)}$,
$|G'/\gamma_{3}(G)|=p$ and hence $|\gamma_{3}(G)|=p^2$. If $|h^G|=p$ for some $h\in G-G'$, then $C_G(h)=\langle h,G'\rangle$ is a maximal abelian subgroup of $G$ and thus $\mathrm{Out}_c(G)=1$; and if $|h^G|=p^3$ for some $h\in G-G'$, then $|C_G(h)|=p^2$ and hence $h\in Z(G)$, which is absurd. We therefore suppose that $|x^G|=p^2$ for all $x\in G-G'$. Then $|\mathrm{Aut}_c(G)|\le p^4$. Let $|\mathrm{Aut}_c(G)|=p^4$. Then, for any $x,y\in G$, there exists an $\alpha\in\mathrm{Aut}_c(G)$ such that $\alpha(a)=x^{-1}ax$ and $\alpha(b)=y^{-1}by$. In particular, there exists an $\alpha\in\mathrm{Aut}_c(G)$ such that $\alpha(a)=b^{-1}ab$ and $\alpha(b)=a^{-1}ba$. Now $(ab)^{-1}\alpha(ab)\in [ab,G]$. But $(ab)^{-1}\alpha(ab)=b^{-1}a^{-1}\alpha(a)\alpha(b)
=b^{-1}wbw^{-1}=[w,b]=v^{-1}$. Thus $v^{-1}\in [ab,G]$. Let $v^{-1}=[ab,h]$, where $h\in G$. Observe that $uv=[a,w][b,w]=[ab,w]$ and therefore $u=(uv)v^{-1}=[ab,w][ab,h]=[ab,wh]\in [ab,G]$. Thus $\gamma_3(G)\le [ab,G]$. But since $|[ab,G]|=|(ab)^G|=p^2$, $\gamma_3(G)=[ab,G]$. Thus $w^{-1}=[b,a]=[ab,a]\in [ab,G]=\gamma_3(G)$. This is a contradiction and hence $\mathrm{Out}_c(G)=1$. This completes the proof. \hfill $\Box$

\begin{thm}
Let $G$ be a finite non-abelian $p$-group of order $p^{5}$. Then  $\mathrm{Out}_{c}(G) \neq 1 $ 
if and only if $|Z(G)|=p,\;Z(G)<G'$
and either $(i)\;cl(G)=3$ and $d(G)=3$ or $(ii)\;cl(G)=4$ and $Z(G) \subseteq [x,G]$ for all $x \in G - G'$.
\end{thm}
{\bf Proof.} First suppose that $\mathrm Out_{c}(G) \neq 1 $. Then $|Z(G)|=p$ by Proposition 2.2 and hence $Z(G) < G'$. It follows that $G$ is purely non-abelian and $cl(G)$ is either 3 or 4. Suppose that $cl(G)=3$. Since $|G'|=p^{2}$ or $p^{3}$,
$|\mathrm{Aut}_z(G)|=|\mathrm{Hom}(G/G',Z(G))|\leq p^{3}$ by \cite[Theorem 1]{AdnYen}. Also, since $G/Z(G)$ is a class 2 group of order $p^{4}$, $|Z(\mathrm{Inn}(G))|=p^{2}$. It now follows by Lemma 2.1 that
$$|\mathrm{Aut}_{c}(G)|= p^{2}|\mathrm{Aut}_c(G) \cap \mathrm{Aut}_{z}(G)| \leq p^{5}.$$
Since $|\mathrm{Aut}_{c}(G)| > |\mathrm{Inn}(G)|=p^{4}$, $|\mathrm{Aut}_z(G)|=p^{3}$ and hence $d(G)=3$. Next suppose that $cl(G)$ = 4. Then $G$ is of maximal class and hence $|G'|=p^{3}$, $|Z(\mathrm{Inn}(G))|=p$, $d(G)=2$ and
$|\mathrm{Aut}_z(G)|=|\mathrm{Hom}(C_{p}\times C_{p} ,C_{p})|= p^{2}$. Thus 
$$|\mathrm{Aut}_{c}(G)|= p^{3}|\mathrm{Aut}_c(G)\cap \mathrm{Aut}_{z}(G)| \leq p^{5},$$ 
by Lemma 2.1. It follows that $|\mathrm{Aut}_c(G)\cap \mathrm{Aut}_{z}(G)|=p^{2}$ and hence $\mathrm{Aut}_z(G) \leq \mathrm{Aut}_{c}(G)$. Let $a\in G - G'$ and let $G=\langle a,b\rangle$. Suppose that all central automorphisms of $G$ fix $a$. Then each $1 \neq \alpha \in \mathrm{Aut}_z(G)$ would have to move $b$. Since $|\mathrm{Aut}_{z}(G)|=p^{2}$, this would require that $|Z(G)|=p^{2}$, which is not so. Thus, there exists a central automorphism $\alpha$ such that $\alpha(a) = az$ for some $1 \neq z\in Z(G)$. On the other hand, since $\alpha$ is class-preserving, $\alpha(a) = g^{-1}ag$
for some $g\in G$. It follows that $z=[a,g]$ and, since $z^{n} = [a,g^{n}]$ for all $n\geq 1$,
$Z(G)\subseteq [x,G]$ for all $x\in G - G'$.

To prove the converse, we first suppose that $|Z(G)|=p$, $Z(G)< G'$, $cl(G)=4$ and $Z(G)\subseteq [x,G]$ for all $x \in G - G'$. Since $G$ is of maximal class, $|Z(\mathrm{Inn}(G))|=p$, $|G'|=p^{3}$, $d(G)=2$ and hence $|\mathrm{Aut}_z(G)|=|\mathrm{Hom}(C_{p}\times C_{p} ,C_{p})|= p^{2}$. Also, since $Z(G)\subseteq [x,G]$ for all $x\in G - G'$, $\mathrm{Aut}_z(G) \leq \mathrm{Aut}_{c}(G)$. It follows by Lemma 2.1 that $| \mathrm{Aut}_{c}(G)|=p^{5}>|\mathrm{Inn}(G)|$. Next suppose that $|Z(G)|=p$, $Z(G)< G'$ and $cl(G)=d(G)=3$. Then $|G'|=|\Phi(G)|=p^{2}$ and hence $|\mathrm{Aut}_z(G)|=|\mathrm{Hom}(G/G',Z(G))|= p^{3}$. It follows from \cite[Theorems 4.7 and 5.1]{TanMor} that $(G,Z(G))$ is a Camina pair. Thus $Z(G)\subseteq [x,G]$ for all $x\in G - Z(G)$ and hence $\mathrm{Aut}_z(G) \leq \mathrm{Aut}_{c}(G)$. Since $|G/Z(G)|=p^{4}$ and $cl(G/Z(G))=2$, $|Z(\mathrm{Inn}(G))|=p^{2}$ and thus $|\mathrm{Aut}_{c}(G)|=p^{5} > |\mathrm{Inn}(G)|$ by Lemma 2.1. This proves the theorem. \hfill $\Box$\\

As a consequence of Theorem 2.3, we can now obtain the following results of Yadav \cite[Section 5]{Yad2008} and Kalra and Gumber \cite[Theorem 4.2]{KalGum}. In \cite{Jam}, the groups of order $p^5$, where $p$ is an odd prime, are divided into ten isoclinism families; and in \cite{HalSen}, the groups of order 32 are divided into eight isoclinism families. The isoclinism families appearing in Corollary 2.4 and Corollary 2.5 are respectively from \cite{Jam} and \cite{HalSen}. The $i$-th family is denoted as $\Phi_i$. We would like to remark here that the derived groups of all the groups in the tenth family of \cite{Jam} are elementary abelian  for $p\ge 5$, but, for $p=3$, the derived groups are not elementary abelian because $|\alpha_2|=9$. It follows that, for $p=3$, the groups $H_1=\Phi_{10}(2111)a_0, H_2= \Phi_{10}(2111)a_1$ and $H=\Phi_{10}(1^5)$ are not in the tenth family.

\begin{cor}
Let $G$ be a finite $p$-group of order $p^5$, where $p$ is an odd prime. Then $\mathrm{Out}_c(G)\ne 1$ if and only if $G$ is isomorphic to one of the groups in $\Phi_7$ or one of the groups in $\Phi_{10}$ for $p\ge 5$ or to $H,H_1$ or $H_2$.
\end{cor}
{\bf Proof.} If $G$ is any group from the first six families, then either $cl(G)<3$ or $|Z(G)|>p$. There are two isoclinism families $\Phi_7$ and $\Phi_8$ consisting of groups of class 3; and there are two isoclinism families $\Phi_9$ and $\Phi_{10}$ consisting of groups of class 4. The only group in eighth family is $\Phi_8(32)$  with $d(\Phi_8(32))=2$. Any group $G$ in the seventh family is generated by $\alpha ,\alpha_1, \beta$, and it is easy to see that the $p$-th power of any of these generators is either 1 or is in $G'$. It thus follows that $|\Phi(G)|=|G'|=p^2$ and consequently $d(G)=3$. Any group $G$ in the ninth family is minimally generated by $\alpha$ and $\alpha_1$ with abelian commutator subgroup $G'=\langle [\alpha_1,\alpha],\gamma_3(G)\rangle=\langle \alpha_2,\gamma_3(G)\rangle$ and $\gamma_3(G)=\langle [\alpha_1,\alpha_2],[\alpha,\alpha_2],Z(G)\rangle=\langle\alpha_3^{-1},Z(G)\rangle$. Thus $\alpha_1$ commutes with $G'$ and hence $|\alpha_{1}^G|=p$. It is easy to see that any element $g$ of $G$ is of the form $g'\alpha_{1}^k\alpha^j$, where $0\le j,k\le p-1$ and $g'\in G'$, and hence $[\alpha_1,g]=[\alpha_1,g'\alpha_{1}^k\alpha^j]=[\alpha_1,\alpha^j]\in G'-\gamma_3(G)$. It follows that $Z(G)$ is not contained in $[x,G]$ for all $x\in G-G'$. Let $G$ be any group in the tenth family. Any element $g\in G-G'$ is of the form $g=g'\alpha_{1}^{l}\alpha^m$, where $0\le l,m\le p-1$ and $g'\in G'$. If $m\ne 0$, then $[g,\alpha_3]=[\alpha^m,\alpha_3]\in Z(G)$. Let $m=0$ and $l\ne 0$. It is easy to see, by induction, that $[\alpha_{1}^{l},\alpha_2]=[\alpha_1,\alpha_2]^{l}z_{1}z_2\ldots z_{l-1}$, where $z_1,z_2,\ldots ,z_{l-1}\in Z(G)$. Then $[g,\alpha_2]=[g'\alpha_{1}^{l},\alpha_2]=[\alpha_{1}^{l},\alpha_2]=[\alpha_1,\alpha_2]^l\in Z(G)$, because $[\alpha_1,\alpha_2]\in Z(G)$. It follows that $Z(G)\subseteq [x,G]$ for all $x\in G-G'$.  \hfill $\Box$\\

There are, in all, fifty one groups of order 32. Sag and Wamsley \cite{SagWam} have given minimal presentations of these groups. As mentioned in \cite{SagWam}, the groups are in the same order in \cite{HalSen} and \cite{SagWam}. We denote the $i$-th group as $G_i$.

\begin{cor}
Let $G$ be a finite group of order $32$. Then $\mathrm{Out}_c(G)\ne 1$ if and only if either $G$ is isomorphic to $G_{44}$ or isomorphic to $G_{45}$.
\end{cor}
{\bf Proof.} First forty three groups are divided into first five families. The families $\Phi_1,\Phi_2,\Phi_4$ and $\Phi_5$ contain groups of class $\le 2$. The third family contains ten groups $G_{23}-G_{32}$ of class 3, and each of these groups has center of order 4. The sixth family contains groups $G_{44}$ and $G_{45}$. Both of these groups are of class 3, rank 3, and with centers of order 2. The seventh family contains groups $G_{46}-G_{48}$ of class 3 and rank 2. The last family contains 2-generated groups $G_{49}-G_{51}$ of maximal class. It is clear from \cite{SagWam} that in each of these groups, one generator, say $x$, is of order 16 and hence $|x^G|=2$. If $y$ is another generator, then $y^{-1}xy$ is respectively $x^{15},x^7$ and $x^{15}$ in the group $G_{49}, G_{50}$ and $G_{51}$. The center of each of these groups is $\{1,x^8\}$. For any element $g$ in these groups, if $[x,g]=x^8$, then $g^{-1}xg=x^9$, which is not so. \hfill $\Box$

\end{document}